\documentstyle[amscd,12pt,psamsfonts]{amsart}

\newtheorem{defn}{Definition}
\newtheorem{thm}[defn]{Theorem}

\newtheorem{cor}[defn]{Corollary}
\newtheorem{lem}[defn]{Lemma}
\newtheorem{prop}[defn]{Proposition}

\theoremstyle{remark}
\newtheorem{rem}{Remark}
\theoremstyle{remark}
\newtheorem{exam}{Example}

\numberwithin{equation}{section}
\numberwithin{defn}{section}

\begin{document}


\newcommand\spk{{\operatorname{Spec}}(k)}
\renewcommand\sp{\operatorname{Spec}}
\renewcommand\sf{\operatorname{Spf}}
\newcommand\proj{\operatorname{Proj}}
\newcommand\aut{\operatorname{Aut}}
\newcommand\grv{{\operatorname{Gr}}(V)}
\newcommand\gr{\operatorname{Gr}}
\newcommand\glv{{\operatorname{Gl}}(V)}
\newcommand\glve{{\widetilde{\operatorname{Gl}}(V)}}
\newcommand\gl{\operatorname{Gl}}
\renewcommand\hom{\operatorname{Hom}}
\renewcommand\det{\operatorname{Det}}
\newcommand\detd{\operatorname{Det}^\ast}
\newcommand\im{\operatorname{Im}}
\newcommand\res{\operatorname{Res}}
\newcommand\limi{\varinjlim}
\newcommand\limil[1]{\underset{#1}\varinjlim\,}
\newcommand\limp{\varprojlim}
\newcommand\limpl[1]{\underset{#1}\varprojlim\,}

\renewcommand\o{{\mathcal O}}
\renewcommand\L{{\mathcal L}}
\newcommand\B{{\mathcal B}}
\renewcommand\c{{{\mathcal C}}^\punto}
\renewcommand\P{{\mathbb P}}
\newcommand\Z{{\mathbb Z}}
\newcommand\A{{\mathbb A}}
\newcommand\M{{\mathcal M}}

\newcommand\w{\widehat}
\renewcommand\tilde{\widetilde}

\newcommand\iso{@>{\sim}>>}
\renewcommand\lim{\limpl{A\in\B}}
\newcommand\fu{\underline}
\newcommand\kz{\fu{k((z))}^*}
\font\peqq=cmr10 at 7pt
\font\peq=cmr10 at 8pt
\font\pun=cmsy10 at 5pt
\newcommand\punto{{\bullet}}
\newcommand\beq{
      \setcounter{equation}{\value{defn}}\addtocounter{defn}1
      \begin{equation}}

\renewcommand{\thesubsection}{\thesection.\Alph{subsection}}

\title [Equations of the moduli of pointed curves]
{Equations of the moduli of pointed curves \\ in the infinite
Grassmannian}

\author[J. M. Mu\~noz \and F. J. Plaza ] {J. M. Mu\~noz
Porras \\ \and \\  F. J. Plaza Mart\'{\i}n\\  \medskip\tiny Departamento
de Matem\'aticas \\ Universidad de Salamanca}

\address{ Departamento de Matem\'aticas \\ Universidad de
Salamanca \\ Plaza de la Merced 1-4\\
Salamanca 37008. Spain.}

\thanks{
This work is partially supported by the
CICYT research contract n. PB96-1305 and Castilla y Le\'on regional
goverment contract SA27/98.}

\email{jmp@@gugu.usal.es}
\email{fplaza@@gugu.usal.es}

\begin{abstract}
The main result of this paper is the explicit computation of the
equations defining the moduli space of triples $(C,p,\phi)$, where $C$ is
an integral and complete algebraic curve, $p$ a smooth rational point
and $\phi$ a certain isomorphism. This is achieved by
introducing algebraically infinite Grassmannians, tau and Baker-Ahkiezer
functions and by proving an Addition Formula for tau functions.
\end{abstract}

\maketitle


\setcounter{tocdepth}1
\tableofcontents

\section{Introduction}

The Krichever morphism gives an immersion of the moduli space,
$\M_\infty$, of triples $(C,p,\phi)$ (where $C$ is an algebraic curve,
$p$ a smooth point of $C$ and $\phi$ a  certain isomorphism) in the
infinite Grassmannian (see \cite{SS,SW}). The aim of this paper is
to give an explicit system of equations defining the subscheme
$\M_\infty$ of $\grv$, this problem is solved in \S6.

We have adopted an algebraic point of view (\cite{AMP}) and most of the
results of this paper are valid over arbitrary base fields.
Consequently, we have included several paragraphs addressing certain
aspects of the theory of soliton equations and infinite Grassmannians
for arbitrary base fields (these facts are well known when the base
field is $\mathbb C$,
\cite{SW,DJKM,KSU}). This allows us to give the foundations for a theory
of soliton equations valid over arbitrary base fields, extending the
previous results of G. Anderson (\cite{An}) for the case of $p$-adic
fields. We hope that the techniques developed in this paper will clarify
the ``arithmetic properties'' of the theory of KP equations.

When the base field is $\mathbb C$, our equations for $\M_\infty$ (as a
subscheme of the infinite Grassmannian, $\grv$) are equivalent to a
system of infinite partial differential equations {\ref{thm:pde-mod}}
which are different from the equations of the KP hierarchy. To clarify the
relation between the KP hierarchy and our  equations {\ref{thm:pde-mod}},
let us consider the chain of closed immersions:

$$\M_\infty\overset{\text{\peqq Krichever}}\hookrightarrow \grv
\overset{\text{\peqq Pl\"ucker}}\hookrightarrow {\mathbb P}^\infty$$
where ${\mathbb P}^\infty$ is a suitable infinite dimensional projective
space. It is well known (\cite{Pl2,SS}) that the image of $\grv$ in
${\mathbb P}^\infty$ is defined by the Pl\"ucker equations, which are
 equivalent to the KP hierarchy ($\operatorname{char}(k)=0$). Then, the
image of
$\M_\infty$ in ${\mathbb P}^\infty$ will be defined by the following
system of differential equations:

$$\left\{\aligned
&\text{ the KP equations (given in Theorem {\ref{thm:KP-eq}}),} \\
&\text{ the p.d.e.'s  given in Theorem {\ref{thm:pde-mod}},} \\
&\text{ the p.d.e.'s given in Corollary {\ref{cor:pde-tau}}.3.}
\endaligned\right.$$

In particular, we deduce a characterization ({\ref{cor:pde-tau}}), in
terms of partial differential equations, of the infinite formal series
$\tau(t)\in {\mathbb C}\{\{t_1,t_2,\dots\}\}$, which are the
$\tau$-functions of a triple $(C,p,\phi)\in\M_\infty$.

Let us note that the results of Shiota (\cite{Sh}) give a necessary
and suficient condition for a theta function of a principally polarized
abelian variety to be the theta function of a Jacobian, but do not solve
the problem of characterizing, in terms of differential equations, the
formal $\tau$-functions defined by Jacobian theta
functions; this problem is solved in Corollary {\ref{cor:pde-tau}}.

The paper is organized as follows. A survey on infinite Grassmannians is
given in \S2. In \S3 the algebraic analogue of the group of maps
$S^1\to{\mathbb C}^*$ of Segal-Wilson is constructed  and  interpreted
as the Jacobian of the formal curve. The action of that group on the
Grassmannian is used in \S4 in order to introduce tau and Baker-Ahkiezer
functions. At the end of this section the Addition Formulae for tau
functions are stated and proved.

Equations defining the infinite Grassmannian of $k((z))$ in a suitable
infinite dimensional projective space are computed in \S5. Over an
arbitrary field it is defined by the well known Bilinear Residue
Identity (equation \ref{eq:residue}) whose proof is an easy consequence
of Theorem \ref{4:thm:BA}. When the base field is ${\mathbb C}$ and
using Schur polynomials, this identity turns out to be equivalent to the
KP Hierarchy.

The last section, \S6, contains the main result of this paper. A
relative generalization of the  Krichever map allows us to obtain a
closed immersion of the moduli functor of pointed curves, whose rational
points are the triples $(C,p,\phi)$, into the infinite
Grassmannian of $k((z))$. We prove then its representability and give a
characterization that permits to compute its equations.

We are very grateful to the referee for his comments, which help us to
improve the paper and to remove some initial mistakes.

\section{Infinite Grassmannians and Determinant Bundles}

\subsection{Infinite Grassmannians}\label{subsect:grass}

In order to define the ``infinite'' Grassmannian of a vector space $V$
(over a field $k$) one should consider some extra structure on it. This
structure consists of a family $\B$ of subspaces of $V$ such that the
following conditions hold:
\begin{enumerate}
\item $A,B\in\B\quad \implies\quad A+B,A\cap B\in\B$,
\item $A,B\in\B\quad \implies\quad \dim(A+B)/A\cap B <\infty$,
\item $\cap_{A\in\B} A=(0)$,
\item the canonical homomorphism $V\to\limpl{A\in\B}V/A$ is an
isomorphism,
\item the canonical homomorphism $\limil{A\in\B}A\to V/B$ is a
surjection (for $B\in\B$).
\end{enumerate}

Let us interpret these conditions. First, $V$ is endowed with a linear
topology where a family of neigbourhoods of $(0)$ is precisely $\B$.
Then, the last three claims mean that the topology is separated; $V$ is
complete and every finite dimensional subspace of $V/B$ is a
neigbourhood of $(0)$.

\begin{exam}\label{2:exam}
In the study of the moduli space of pointed curves the fundamental
example is $V=k((z))$ and $\B$ consists of the set of subspaces
$A\subseteq V$ containing $z^n\cdot k[[z]]$ as a subspace of finite
codimension (for an integer $n$).
\end{exam}

\begin{exam}
Other examples of pairs $(V,\B)$ satisfying the above requeriments are
($V$ an arbitrary $k$-vector space):
\begin{itemize}
\item $(V, \B:=\{(0)\})$;
\item $V$ and $\B$ the set of all finite dimensional subspaces of $V$.
\end{itemize}
\end{exam}

In the sequel, we shall fix a pair $(V,\B)$ satisfying the above
conditions and a subspace $V_+\in\B$.

Following \cite{AMP} (see also \cite{KSU}), there exists a Grassmannian
scheme
$\gr(V,\B)$, whose rational points are the set:
$$\left\{\begin{gathered}
\text{ subspaces $L\subseteq V$, such that $L\cap V_+$ }\\
\text{ and ${V}/{L+V_+}$ are of finite dimension }
\end{gathered}\right\}$$
(which we shall call discrete subspaces of $V$)  that coincides with the
usual infinite Grassmannian defined by Pressley-Segal in \cite{PS},
and Segal-Wilson in \cite{SW}.

In order to construct this scheme, we shall give its functor of points
$\fu{\gr}(V,\B)$ and prove that it is representable in the category of
$k$-schemes. To this end, we need some notation: given a
morphism $T\to S$ of $k$-schemes a sub-$\o_S$-module $B\subseteq
V_S:=V\otimes_k\o_S$, we denote:
\begin{itemize}
\item $\w B_T:=\lim (B/(A_S\cap B) \underset k\otimes\o_T)$.
\item $\w{(V/B)}_T:=\lim ((V_S/(A_S+B))\underset k\otimes\o_T)$.
\end{itemize}

\begin{defn}\label{defn:discrete}
Given a $k$-scheme $S$, a discrete
submodule of $\w V_S$ is a sheaf of quasi-coherent
$\o_S$-submodules $\L \subset \w V_S$ such that: $\L_T\subset
\w V_T$ for every morphism $T\to S$; and, for each
$s\in S$ there exists an open neighborhood
$U$ of $s$ and a commensurable $k$-vector subspace $B\in\B$ such
that $\L_U\cap \w B_U$ is free of finite type and
$\w V_U/\L_U+ \w B_U=0$.
\end{defn}

\begin{defn}
The Grassmannian functor of a pair $(V,\B)$, $\fu{\gr}(V,\B)$, is the
contravariant functor over the category of
$k$-schemes defined by:
$$\fu{\gr}(V,\B)(S)=\left\{
 \text{discrete sub-$\o_S$-modules of $\w V_S$}\right\}$$
\end{defn}

\begin{rem}
Note that if $V$ is a finite dimensional $k$-vector space and
$\B=\{(0)\}$, then $\fu{\gr}(V,\{(0)\})$ is the usual Grassmannian
functor defined by Grothendieck \cite{EGA}~I.9.7.
\end{rem}

\begin{thm}
The functor $\fu{\gr}(V,\B)$ is representable by a reduced and
separated $k$-scheme $\gr(V,\B)$. The discrete sub\-module  corresponding to
the identity:
$$Id\in \fu{\gr} (V,\B)\left(\gr(V,\B)\right)$$
will be called the {\bf universal submodule} and will be denoted by:
$$\L_V\subset \w V_{\gr(V,\B) }$$
\end{thm}

\begin{pf}
The proof is modeled on the Grothendieck construction
of finite Grassmannians \cite{EGA}. Given a vector subspace
$A\in \B $, define the functor $\fu{F_A}$ over the category of
$k$-schemes by:
$$\fu{F_A}(S)=\{\text{ sub-$\o_S$-modules $\L\subset \w V_S$
such that $\L\oplus \w A_S=\w V_S$ }\}$$
and show that it is representable by an affine and integral $k$-scheme,
$F_A$. The properties on $\B$ imply that
$\{\fu{F_A}, A\in \B \}$ is a covering of $\fu{\gr}(V,\B)$ by open
subfunctors (see \cite{AMP}) and the result follows.
\end{pf}

\begin{rem}
In this subsection infinite-dimensional
Grassmannian sche\-mes have been constructed in an abstract way. Choosing
particular vector spaces $(V,\B)$ we obtain different classes of
Grassmannians. Two examples are relevant:
\begin{itemize}
\item $V=k((z))$, $V_+=k[[z]]$. In this case, $\gr(k((z)),k[[z]])$ is
the algebraic version of the Grassmannian constructed by Pressley-Segal,
and Segal-Wilson (\cite{PS,SW}). This Grassmannian is
particularly well adapted for studying problems related to the moduli
space of pointed curves (over arbitrary fields), to the moduli space of
vector bundles and to the KP-hierarchy.
\item Let $(X,\o_X)$ be a smooth, proper and irreducible curve over the
field $k$ and let $V$ be the adeles ring  over the curve and
$V_+=\underset p \prod\w{\o_p}$ (recall the first example). In this
case  $\gr(V,\B)$ is an adelic Grassmannian which will be useful for
studying arithmetic problems over the curve $X$ or problems related to
the classification of vector bundles over a curve (non abelian theta
functions...).
\newline
Instead of adelic Grassmannians, we could define Grassmannians
associated with a fixed divisor on $X$ in an analogous way.
\newline
These adelic Grassmannians  will also be  of interest in the study of
conformal field theories over Riemann surfaces in the sense of
Witten (\cite{W}).
\end{itemize}
\end{rem}


\subsection{Determinant Bundles}\label{subsect:det}

In this subsection we recall from \cite{AMP} the construction the
determinant bundle over the Grassmannian, in the sense  of Knudsen and
Mumford \cite{KM}. This allow us to define determinants
algebraically and over arbitrary fields (for example for $k={\mathbb Q}$
or $k={\mathbb F}_q$).

We shall denote the Grassmannian $\gr(V,\B)$ simply by $\grv$.

\begin{defn}
For each $A\in \B $ and each $L\in\grv(S)$ we define a complex,
$\c_A(L)$, of $\o_S$-modules by:
$$\c_A(L) \equiv \ldots\to 0\to L\oplus \hat A_S
@>\delta>> \hat V_S\to 0\to\ldots$$
$\delta$ being the addition homomorphism.
\end{defn}

It is not difficult to prove that $\c_A(L)$ is a perfect complex of
$\o_S$-modules, and therefore its determinant and its index are
well defined. Recall that the index of a point $L\in\grv(S)$ is the
locally constant function $i_L\colon S\to\Z$ defined by:
$$i_L(s)=\text{  Euler-Poincar\'e characteristic of }
\c_{V_+}(L)\otimes k(s)$$
$k(s)$ being the residual field of the point $s\in S$ (see \cite{KM} for
details).

By easy calculation, we have that if $V$ is a finite-dimensional
$k$-vector space, $\B=\{(0)\}$ and $L\in\grv(S)$, then
$i_L=\operatorname{codim}_k(L)$. And for the general case, the index of
any rational point $L\in\grv(\spk)$ is exactly:
$$\dim_k(L\cap V_+)-\dim_k V/(L+ V_+)$$

Let $\gr^n(V)$ be the subset over which the index takes values
equal to $n\in\Z$. Then, the decomposition of $\grv$ in connected
components is:
$$\grv=\underset{n\in\Z}\coprod\gr^n(V)$$

Given a point $L\in\grv(S)$ and $A\in \B $, we denote by
$\det\c_A(L)$ the determinant sheaf of the perfect complex
$\c_A(L)$ in the sense of \cite{KM}. Note that this determinant does not
depend on $A$ (up to isomorphisms), since for $A,B\in\B$ there exists
a canonical isomorphism:
\beq
\detd\c_A\iso\detd\c_B\otimes\wedge^{max}( A/{A\cap
B})\otimes \wedge^{max}(B/{A\cap B})^*
\label{2:eqn:det-iso}\end{equation}
Hence we define  the determinant bundle over
$\gr^0(V)$,
$\det_V$, as the invertible sheaf:
$$\det\c_{V_+}(\L_V)$$
($\L_V$ being the universal submodule over $\gr^0(V)$).

Observe that the choice of other subspace of $\B$, $A$, instead of $V_+$
might shift the labelling of the connected components by a constant and
modify the determinant bundle by an isomorphism.

Now, if $L\in\gr^0(V)$ is a rational point, and  $L\cap V_+$ and
${V}/{L+ V_+}$ are therefore
$k$-vector spaces  of the same dimension, then one has an
isomorphism:
$$\det_V(L)\simeq
\wedge^{max}(L\cap V_+)\otimes\wedge^{max}({V}/{(L+V_+)})^\ast$$
That is, our determinant coincides, over the rational points, with
the determinant bundles of Pressley-Segal and Segal-Wilson (\cite{PS,SW}).
Furthermore, this construction gives the usual determinant bundle when
$V$ is finite-dimensional.

\subsection{Sections of the Determinant bundle}

It is well known that the determinant bundle has no
global sections. We shall therefore explicitly
construct global sections of the dual of the determinant
bundle over the connected component $\gr^0(V)$ of index
zero.

We use the following notations: $\wedge^\punto E$ is the
exterior algebra of a $k$-vector space  $E$, $\wedge^r E$
its component of degree $r$, and $\wedge E$ is the component
of higher degree when $E$ is finite-dimensional.
Given a perfect complex $\c$ over  $k$-scheme $X$, we shall
write $\detd\c$ to denote the dual of the invertible
sheaf $\det\c$.

To explain how global sections of the invertible sheaf $\detd\c$ can be
constructed, recall that if $f:E\to F$ is a
homomorphism between finite-dimensional $k$-vector spaces of equal
dimension, then it induces a homomorphism:
$$\wedge(f):k\to \wedge F\otimes (\wedge E)^*$$
Thus, considering $E @>{f}>> F$ as a perfect
complex, $\c$, over $\spk$, we have defined a {\bf
canonical section} $\wedge(f)\in H^0(\spk,\detd\c)$.

Let us now consider a perfect complex
$\c\equiv(E@>{f}>>F)$ of sheaves of $\o_X$-modules over
a $k$-scheme $X$, with Euler-Poincar\'e characteristic
${\mathcal X}(\c)=0$. Using the above argument, we
construct a canonical section $det(f\vert_U)\in H^0(U,\detd\c)$ for
every open subscheme of $X$, $U$, over which $\c$ is quasi-isomorphic to
a complex of finitely-generated free modules.

Since the construction is canonical, the functions:
$$g_{UV}:=
det(f\vert_U)\vert_{U\cap V}\cdot det(f\vert_V)\vert_{U\cap V}^{-1}$$
satisfy the cocycle condition. Then, these local sections glue
$\{det(f\vert_U)\}$ and give a canonical global section:
$$\det(f)\in H^0(X,\detd\c)$$
(see \cite{AMP} for more details). If the
complex
$\c$ is acyclic, one has an isomorphism:
$$\aligned
\o_X &\iso\det^\ast\c\\
1 &\mapsto det(f)
\endaligned$$

Let us consider the perfect complex $\c_A\equiv(\L\oplus  A
@>{\delta_A}>>  V)$ over $\grv$ defined in {\ref{subsect:det}} ($\L$
being the universal discrete submodule over $\grv$) for a given
$A\in\B$. Since $\c_A\vert_{F_A}$ is acyclic, one then has an
isomorphism:
$$\aligned
\o_X\vert_{F_A}& @>\sim>>
\detd\c_A\vert_{F_A}\\
1 &\mapsto s_A=det(\delta_A)\vert_{F_A}
\endaligned$$

By the above argument it is easy to prove that the section
$s_A\in H^0(F_A,\detd\c_A)$ can be extended in a
canonical way to a global section of $\detd\c_A$ over $\gr^0(V)$, which
 will be called the {\bf canonical section $\omega_A$ of
$\detd\c_A$}. (We restrict ourselves to cases where
$F_A\subseteq \gr^0(V)$, or, what amounts to the same, $\dim_k( A/{A\cap
V_+})-\dim_k({V_+}/{A\cap V_+})=0$).

This result allows us to compute many global sections
of $\detd_V=\det\c_{V_+}$ over $\gr^0(V)$: given $A\in \B $ such that
$F_A\subseteq\gr^0(V)$, the isomorphism $\detd\c_A\iso\detd_V$  is not
canonical (recall the formula \ref{2:eqn:det-iso}).
Therefore, the determination of an isomorphism $\detd\c_A\iso\detd_V$
depends on the choice of bases for the vector spaces
$ A/{A\cap V_+}$ and ${V_+}/{A\cap V_+}$.
For a detailed discussion of the construction of sections see
\cite{AMP,Pl2}).

\subsection{Computations for the infinite
Grassmannian: $V=k((z))$}\label{subsect:comp-inf-grass}

Let $(V,\B)$ be as in Example \ref{2:exam} and take $V_+=k[[z]]$.

Let $\mathcal S$ be the set of Young's diagrams (also called Maya or
Ferrers diagrams) of virtual cardinal zero; or equivalently,
the sequences
$\{s_0,s_1,\ldots\}$ of integer numbers satisfying the following
conditions:
\begin{enumerate}
\item the sequence is strictly increasing,
\item  there exists $s\in\Z$ such that
$\{s,s+1,s+2,\ldots\}\subseteq\{s_0,s_1,\ldots\}$,
\item $\#(\{s_0,s_1,\ldots\}-\{0,1,\ldots\})=
\#(\{0,1,\ldots\}-\{s_0,s_1,\ldots\})$
\end{enumerate}

For notation's sake, we define $e_i:=z^i$. For each $S\in{\mathcal S}$,
let $A_S$ be the vector subspace of $V$ generated by $\{e_{s_i}, i\ge
0\}$. By the third condition one has:
$$\dim_k({A_S}/{A_S\cap V_+})=
\dim_k({V_+}/{A_S\cap V_+})$$
and hence $A_S\in \B $ and $F_{A_S}\subseteq\gr^0(V)$.
Further, $\{F_{A_S},S\in{\mathcal S}\}$ is a covering of
$\gr^0(V)$.

Define linear forms $\{e^*_i\}$ of $V^*$ by the following condition:
$$e^*_i(e_j)\,:=\,\delta_{ij}$$
For each finite set of increasing integers,
$J=\{j_1,\dots,j_r\}$, let us define
$e_J:=e_{j_1}\wedge\ldots\wedge e_{j_r}$ and
$e^*_J:=e^*_{j_1}\wedge\ldots\wedge e^*_{j_r}$.

Given $S\in{\mathcal S}$, choose $J,K\subseteq\Z$ such that
$\{e_j\}_{j\in J}$ is a basis of ${V_+}/{A_S\cap V_+}$ and
$\{e^*_k\}_{k\in K}$ of $({A_S}/{A_S\cap V_+})^*$ . We have seen that
tensor by $e_J\otimes e^*_K$ defines an isomorphism:
$$ H^0(\gr^0(V),\detd\c_{A_S})
@>{\,\otimes(e_J\otimes e^*_K)\,}>>
 H^0(\gr^0(V),\detd_V)$$

\begin{defn}\label{defn:global-section}
For each $S\in{\mathcal S}$, $\Omega_S$ is the global section
of $\detd_V$ defined by:
$$\Omega_S=\omega_{A_S}\otimes e_J\otimes e^*_K$$
We shall denote by $\Omega_+$ the canonical section of
$\detd_V$.
\end{defn}

Let $\Omega({\mathcal S})$ be the $k$-vector subspace of
$ H^0(\gr^0(V),\detd_V)$ generated by the global sections
$\{\Omega_S,S\in{\mathcal S}\}$.

We define the Pl\"ucker morphism:
$$\aligned{\frak p}_V: \gr^0(V) &\to \P\Omega(S)^*:=\proj
\operatorname{Sym}\Omega(S)
\\ L &\mapsto \{\Omega_S(L)\}\endaligned$$
as the morphism of schemes associated to the surjective sheaf
homomorphism:
$$\Omega(S)_{\gr^0(V)}\to\detd_V$$
by the universal property of $\P$.

\begin{rem}\label{rem:grass=plucker}
Once the Pl\"ucker morphism is introduced it can be proved that the
Pl\"ucker equations are in fact the
defining equations for $\gr^0(V)$ when $\operatorname{char}(k)=0$
(\cite{Pl2}).

This is the property that Sato used in \cite{SS} to define his Universal
Grassmann Manifold (UGM); that is, a point of the UGM is a point of an
infinite dimensional projective space (with countable many coordinates)
satisfying all the Pl\"ucker relations.
\end{rem}

\begin{rem}\label{rem:indexn}
For studing the index $n$ connected component of the
Grassmannian,
$\gr^n(V)$ the same constructions are applied since the
homothety of $k((z))$ defined by $z^{-n}$ induces isomorphisms:
$$\begin{gathered}
\gr^n(V)\,\iso\,\gr^0(V) \\
H^0(\gr^0(V),\detd_V)\,\iso \, H^0(\gr^n(V),\detd_n)
\end{gathered}$$
where $\detd_n$ is defined by $\detd\c_{z^n\cdot V_+}$. Let us
distinguish the global section given by the image of $\Omega_+$ and
denote it by $\Omega_+^{n}$.
\end{rem}


\subsection{The Grassmannian of the dual space}\label{subsect:grass-dual}

Let $(V,\B)$ be as usual. Consider $V$ as a linear topological space.
Now, a submodule $L\subseteq \w V_S$ ($S$ a $k$-scheme) carries a linear
topology: that given by $\{L\cap \w A_S\}_{A\in\B}$ as neigbourhoods of
$(0)$. We introduce the following notation:
$$\begin{aligned}
L^* &\,:=\, \hom_{\o_S}(L,\o_S) \\
L^c &\,:=\, \{f\in L^* \text{ continuous }\} \\
L^\circ &\,:=\, \{f\in (\w V_S)^* \,\vert\, f\vert_L\equiv 0 \} \\
L^\diamond &\,:=\, \{f\in (\w V_S)^c \,\vert\, f\vert_L\equiv 0 \}
\end{aligned}$$
where $\o_S$ has the discrete topology.

Observe that given two subspaces $A,B\in\B$ such that $B\subseteq A$ the
following claims hold:
\begin{itemize}
\item there is a canonical isomorphism $A^\diamond / B^\diamond \iso
(A/B)^*$. (This implies that $(A^\diamond + B^\diamond)/A^\diamond \cap
B^\diamond$ is finite dimensional).
\item $(A+B)^\diamond= A^\diamond \cap B^\diamond$.
\item $(A\cap B)^\diamond= A^\diamond + B^\diamond$.
\end{itemize}

Consider the following family of subspaces of $V^c$:
$$\B^\diamond\,:=\,\{A^\diamond\text{ where }A\in\B\}$$

In order to make explicit the meaning of the expression
``Grassmannian of the dual space'', $(V^c,\B^\diamond)$, we need the
following:

\begin{lem}[\cite{Pl2}]\hfill
\begin{enumerate}
\item $V^c=\limil{A\in\B}A^\diamond$;
\item $V=\limil{A\in\B} A$;
\item $\cap_{A\in\B} A^\diamond =(0)$;
\item $V^c=\limpl{A\in\B} V^c/A^\diamond$.
\end{enumerate}
\end{lem}

The Lemma and these considerations imply the following:

\begin{thm}
The family $\B^\diamond$ satisfies the conditions of
\ref{subsect:grass}, and therefore the infinite Grassmannian of the pair
$(V^c,\B^\diamond)$ exists.
\end{thm}

If a subspace $V_+\in\B$ is chosen, then we shall consider the subspace
$V_+^\diamond$ for the pair $(V^c,\B^\diamond)$.

Now, we shall construct a canonical isomorphism between the
Grassmannian of
$V$ and that of $V^c$. The expression of this isomorphism for the
rational points will be that given by incidence:
$$\aligned I:\gr(V,\B) &\longrightarrow \gr(V^c,\B^\diamond) \\
L\quad &\longmapsto\quad L^\diamond \endaligned $$

Let $\L$ be the universal sheaf of $\grv$. Consider the following
sub-$\o_{\grv}$-module of $\w V^c_{\grv}$ defined as:
$$\L^\diamond\,=\,
\{ \omega\in \w V^c_{\grv} \text{ such that }\omega\vert_{\L}\equiv 0
\}$$
Let us check that $\L^\diamond$ is in fact a $\grv$-valued point of
$\gr(V^c)$, and that it does induce the desired morphism. By the
definition of the Grassmannian, this can be done locally. Recall that
$\{F_A\}_{A\in \B }$  is a covering of
$\grv$. Let $\L^\diamond\vert_{F_A}$ be the restriction of
$\L^\diamond$ (as sub-$\o_{\grv}$-module of $\w V^c_{\grv}$) to the open
subscheme $F_A\hookrightarrow \grv$. Since $\w V_{F_A}\simeq
\L\vert_{F_A}\oplus\w A_{F_A}$ canonically, one has:
$$\L^\diamond\vert_{F_A} \,\simeq\, (\w A_{F_A})^c$$
in a canonical way. Recalling the Definition {\ref{defn:discrete}}, the
conclusion follows.

To finish, we compute $I^*\det_{V^c}$. Observe that there exists a
canonical morphism of complexes of $\o_{\grv}$-modules (written
vertically):
$$\CD
\L_V^\diamond\oplus V_+^\diamond @>>> V^*  \\
@VVV @VVV \\
V^c @>>> \L_V^*\oplus V_+^* \endCD$$
(where $\L_V$ is the universal submodule of $\grv$) and one
easily checks  that this is in fact a quasi-isomorphism. Since the
inverse image of the universal submodule of $\gr(V^c)$ is
$\L_V^\diamond$, one has the following formulae:
$$I^*\det_{V^*}\simeq \detd_V\qquad\qquad I^*(i)=-i$$
 ($i$ is the index function).


\section{``Formal Geometry'' of Local Curves}

\subsection{Formal Groups}

We are first concerned with the
algebraic analogue of the group $\Gamma$ (\cite{SW}~\S~2.3) of
continuous maps $S^1\to{\mathbb C}^*$ acting as multiplication operators
over the Grassmannian. The main difference between our definition of the
group $\Gamma$ and the definitions offered in the literature
(\cite{SW,PS}) is that in the algebro-geometric setting the
elements $\sum_{-\infty}^{+\infty}g_k\,z^k$ with
infinite positive and negative coefficients do not make
sense as multiplication operators over $k((z))$. In this sense, our
approach is close to that of \cite{KSU}.

The main idea for defining the algebraic analogue of the
group $\Gamma$ is to construct a (formal) scheme whose set of
rational points is precisely the multiplicative group
$k((z))^*$ (see \cite{AMP}).

\begin{defn}
The contravariant functor, $\kz$, over the category of
$k$-schemes with values in the category of commutative
groups is defined by:
$$S\rightsquigarrow \kz(S)\,:=\, H^0(S,\o_S)((z))^*$$
Where for a $k$-algebra $A$, $A((z))^*$ is the group of
invertible elements of the ring $A((z)):=A[[z]][z^{-1}]$ of Laurent
series with coefficients in $A$.
\end{defn}

Note that for each $k$-scheme $S$ and $f\in\kz(S)$, the function:
$$\aligned S&\to \Z \\
s &\mapsto v_s(f):= \text{order of $f_s\in k(s)((z))$}
\endaligned$$
is locally constant. From this fact we deduce that given an irreducible
affine $k$-scheme, $S=\sp(A)$, we have that:
$$\aligned
\kz(S)=&\coprod_{n\in \Z}\left\{f\in A((z))^* \,\vert\, v(f)=n\right\} \\
=&\coprod_{n\in \Z}\left\{\gathered
\text{series }\,a_{n-r}\,z^{n-r}+\dots+a_n\,z^n+\dots\text{ such
that}\\ a_{n-r},\dots,a_{n-1}\text{ are nilpotent and }a_n\in A^*
\endgathered\right\}
\endaligned$$

\begin{thm}
The subfunctor  $\kz_{red}$ of $\kz$ defined by:
$$S\rightsquigarrow \kz_{red}(S)\,:=\,
\coprod_{n\in \Z}\left\{z^n+\sum_{i> n} a_i\,z^i \quad
a_i\in H^0(S,\o_S)\right\}$$
is representable by a group $k$-scheme whose connected component
of the origin will be denoted by $\Gamma_+$.
\end{thm}

\begin{pf}
It suffices to observe that the functor:
$$S\rightsquigarrow \left\{z^n+\sum_{i> n} a_i\,z^i \quad
a_i\in H^0(S,\o_S)\right\}$$
is representable by the scheme $\sp(\limil{l}
k[x_1,\dots,x_l])=\limpl{l}\A^l_k$, where the group law is given
by the multiplication of series; that is:
$$\aligned
k[x_1,\dots]&\to k[x_1,\dots]\otimes_k k[x_1,\dots] \\
x_i &\mapsto x_i\otimes 1+\sum_{j+k=i}x_j\otimes x_k+1\otimes x_i
\endaligned$$
\end{pf}

\begin{thm}
Let $\kz_{nil}$ be the subfunctor of $\kz$ defined by:
$$S\rightsquigarrow \kz_{nil}(S)\,:=\,
\coprod_{n>0}\left\{\gathered
\text{finite series }\,a_n\,z^{-n}+\dots+a_1\,z^{-1}+1\\
\text{such that $a_i\in H^0(S,\o_S)$ are nilpotent}
\endgathered\right\}$$
There exists a formal group $k$-scheme $\Gamma_-$ such that:
$$\hom_{\text{for-sch}}(S,\Gamma_-)=\kz_{nil}(S)$$
for every $k$-scheme $S$.
\end{thm}

\begin{pf}
Note that $\Gamma_-$ is the direct limit in the
category of formal schemes (\cite{EGA} {\bf I}.10.6.3) of the schemes
 representing the subfunctors:
{\small $$S\rightsquigarrow \Gamma^n_-(S)=\left\{\gathered
\,a_n\,z^{-n}+\dots+a_1\,z^{-1}+1\text{ such that }a_i\in H^0(S,\o_S)\\
\text{ and the ${\text{n}}^{\text{th}}$ power of the ideal
$(a_1,\dots,a_n)$ is zero}
\endgathered\right\}$$}
And its associated ring (that is, as ringed space) is:
$$k\{\{x_1,\dots\}\}=\underset n\limp k[[x_1,\dots,x_n]]$$
the morphisms of the projective system being:
$$\aligned
k[[x_1,\dots,x_{n+1}]]&\to k[[x_1,\dots,x_n]] \\
x_i&\mapsto x_i\qquad\text{for }i=1,\dots,n-1 \\
x_{n+1}&\mapsto 0
\endaligned$$
It is now easy to show that $\Gamma_-=\sf(k\{\{x_1,\dots\}\})$
(with group law given by multiplication of series) satisfies the desired
condition. We shall call the ring $k\{\{x_1,\dots\}\}$ the ring
of ``infinite'' formal series in infinite variables (which is different from
the ring of formal series in infinite variables. For instance,
$x_1+x_2+\ldots\in k\{\{x_1,\dots\}\}$).
\end{pf}

\begin{rem}
Our group scheme $\Gamma$ is the algebraic analogue of the
$\Gamma$ group of Segal-Wilson \cite{SW}. Note that the
indices ``-'' and ``+'' do not coincide with the Segal-Wilson
notations. Replacing $k((z))$ by $k((z^{-1}))$, we obtain their
notation.
\end{rem}

Let us define the exponential maps for the groups $\Gamma_-$ and
$\Gamma_+$. Let $\A_n$ be the $n$ dimensional affine space over
$\spk$ with the additive group law, and ${\hat\A}_n$
the formal group obtained as the completion of $\A_n$ at
the origin. We define ${\hat\A}_\infty$ as the formal group scheme
$\limil{n}{\hat\A}_n$. Obviously it holds that:
$${\hat\A}_\infty=\sf k\{\{y_1,\dots\}\}$$
with the additive group law.

\begin{defn}
If the characteristic of $k$ is zero, the exponential map for
$\Gamma_-$ is the following isomorphism of formal group schemes:
$$\aligned
{\hat\A}_\infty & @>{\exp}>> \Gamma_- \\
\{a_i\}_{i>0} &\mapsto \exp(\sum_{i>0}a_i\,z^{-i})
\endaligned$$
This is the morphism induced by the ring homomorphism:
$$\aligned
 k\{\{x_1,\dots\}\}& @>{\qquad\exp^*\qquad}>>
k\{\{y_1,\dots\}\}\\  x_i &\mapsto \text{ coefficient of $z^{-i}$
in the series }
\exp(\sum_{j>0}y_j\,z^{-j})\endaligned$$
\end{defn}

\begin{defn}\label{defn:exp-gamma-minus}
If the characteristic of $k$ is positive, the exponential map for
$\Gamma_-$ is the following isomorphism of formal schemes:
$$\aligned {\hat\A}_\infty &\to \Gamma_- \\
\{a_i\}_{i>0}&\mapsto \prod_{i>0}(1-a_i\,z^{-i})\endaligned$$
which is the morphism induced by the ring homomorphism:
$$\aligned k\{\{x_1,\dots\}\}&@>\exp^*>> k\{\{y_1,\dots\}\}\\
x_i &\mapsto \text{ coefficient of $z^{-i}$ in the series }
\prod_{i>0}(1-a_i\,z^{-i})\endaligned$$
\end{defn}

Note that this latter exponential map is not a
isomorphism of groups. Considering over ${\hat\A}_\infty$ the
law group induced by the isomorphism, $\exp$, of formal schemes,
we obtain the Witt formal group law.

The exponential map for $\Gamma_+$ is defined in a analogous way; one
only has to replace $z^{-i}$ by $z^i$ in the above expressions. (See
\cite{B} for the connection of these definitions and the
Cartier-Dieudonn\'e theory). The following property gives the structure
of $\kz$.

\begin{thm}
The natural morphism of functors of groups over the category of
$k$-schemes:
$$\fu{\Gamma_-}\times\fu{{\mathbb G}_m}\times\fu{\Gamma_+}\to\kz$$
is injective and for $char(k)=0$ gives an isomorphism with
$\kz_0$ (the connected component of the origin in the
functor of groups $\kz$). The functor on groups $\kz_0$ is therefore
``representable'' by the (formal) $k$-scheme:
$$\Gamma=\Gamma_-\times{\mathbb G}_m\times\Gamma_+$$
\end{thm}

\subsection{``Formal Geometry''}

It should be noted that the formal group scheme $\Gamma_-$ has
properties formally analogous to the Jacobians of algebraic
curves: one can define formal Abel maps and prove formal
analogues of the Albanese property of the Jacobians of smooth
curves.

Let $\hat C=\sf(k[[t]])$ be a formal curve. The Abel
morphism of degree $1$ is defined as the morphism of formal schemes:
$$\phi: \hat C\to \Gamma_-$$
given by the $\hat C$-valued point of $\Gamma_-$:
$\phi(t)=(1-\frac{t}z)^{-1}=1+\sum_{i>0}^{}\frac{t^i}{z^i}$;
that is, the morphism induced by the ring homomorphism:
$$\aligned k\{\{x_1,\dots\}\}&\to k[[t]]\\
x_i\,&\mapsto t^i\endaligned$$

Note that the Abel morphism is the algebro-geometric
version of the function $q_\xi(z)$ used by Segal-Wilson (\cite{SW}
page~32) to study the Baker-Akhiezer function.

Let us explain further why we call $\phi$ the
``Abel morphism'' of degree 1. If $char(k)=0$, composing
$\phi$ with the inverse of the exponential map, affords:
$$\bar\phi:\hat C @>{\phi}>> \Gamma_- @>{\exp^{-1}}>>
{\hat\A}_\infty$$
and since
$(1-\frac{t}z)^{-1}=\exp(\sum_{i>0}^{}\frac{t^i}{i\,z^i})$ (see
\cite{SW} page~33), $\bar\phi$ is the morphism defined by the ring
homomorphism:
$$\aligned
k\{\{y_1,\dots\}\} &\to k[[t]]\\
y_i&\mapsto \frac{t^i}i
\endaligned$$
or in terms of the functor of points:
$$\aligned\hat C& @>{\bar\phi}>> {\hat\A}_\infty\\
t&\mapsto \{t,\frac{t^2}2,\frac{t^3}3,\dots\}
\endaligned$$
Observe that given the basis $\omega_i=t^i\,dt$ of the
differentials $\Omega_{\hat C}=k[[t]]dt$, $\bar\phi$ can be
interpreted as the morphism defined by the ``abelian integrals''
over the formal curve:
$$\bar\phi(t)=\left(
\int_0^t\omega_0,\int_0^t\omega_1,\dots,\int_0^t\omega_i,\dots
\right)$$
which coincides precisely with the local equations of the Abel
morphism for smooth algebraic curves over the field of complex
numbers.

\begin{rem}\label{rem:uni-ele}
The above introduced notation is also motivated by the following two facts:
\begin{itemize}
\item  $(\Gamma_-,\phi)$ satisfies the Albanese property for $\hat
C$; that is, every morphism
$\psi:\hat C\to X$ in a commutative group scheme (which sends the unique
rational point of $\hat C$ to the $0\in X$) factors through the
Abel morphism and a homomorphism of groups $\Gamma_-\to X$. This
property follows from the following fact: the direct limit
$\limi_{n}S^n\w C$ exists (as a formal scheme) and it is naturally
isomorphic to $\Gamma_-$).
\item  Observe that for each element:
$$u\in\Gamma_-(S)\subseteq\fu{k((z))^*}(S)=H^0(S,\o_S)((z))^*$$
 we can define a fractionary ideal of the formal curve $\hat C_S$ by:
$I_u=u\cdot\o_S((z))$. We can therefore interpret the
formal group $\Gamma_-$ as a kind of Picard scheme over the formal
curve $\hat C$ (see \cite{C,B2}). The universal element of
$\Gamma_-$ is the invertible element of
$\kz(\Gamma_-)$ given by:
$$v(x,z)=1+\underset {i \geq 1}\sum x_i\,z^{-i}\in k((z))\hat \otimes
k\{\{x_{ 1},x_{ 2},\dots\}\}$$
This universal element will be the formal analogue of the universal
invertible sheaf for the formal curve $\hat C$.
\end{itemize}
\end{rem}

\section{$\tau$-functions and Baker-Akhiezer functions}

The first part of this section is devoted to  defining the
$\tau$-function and the Baker-Akhiezer function algebraically over  an
arbitrary base field $k$. We then prove an analogue for $\tau$ of the
Addition
 Formulae for theta functions that will allow a characterization of the
Baker-Akhiezer function,  which is quite close to proposition 5.1 of
\cite{SW}.

Following on with the analogy between the groups $\Gamma$ and
$\Gamma_-$ and the Jacobian of the smooth algebraic curves, we shall
perform the well known constructions for the Jacobians of the algebraic
curves for the formal curve $\hat C$ and the group $\Gamma$: Poincar\'e
bundle over the dual Jacobian and the universal line bundle over the
Jacobian. In the formal case, these constructions are essentially
equivalent to defining the
$\tau$-functions and the Baker-Akhiezer functions.

Let us denote by $\gr^0(V)$ index zero connected component of the
Grassmannian of
$V=k((z))$  and by
$\Gamma$ the group  $\Gamma_-\times {\mathbb G}_m \times \Gamma_+$. Let
$$\Gamma \times \gr^0(V) \overset \mu \longrightarrow \gr^0(V) $$
be the action of $\Gamma$ over the Grassmannian by homotheties. We define
the Poincar\'e bundle over $\Gamma\times \gr^0(V)$ as the invertible sheaf:
$${\mathcal P}=\mu^*\det_V^*$$
$p_2:\Gamma\times\gr^0(V) \longrightarrow \gr^0(V)$ being the natural
projection.

For each point $U\in \gr^0(V)$, let us define the Poincar\'e bundle over
$\Gamma\times \Gamma$ associated with $U$ by:
$${\mathcal P}_U=(1\times \mu_U)^*{\mathcal
P}=m^*({\mu_U}^*\det_V^*)$$ where $\mu_U:\Gamma\to  \Gamma(U)
\subset \gr^0(V)$ is the action of
$\Gamma$ on the orbit of $U$ and $ m :\Gamma\times\Gamma  \to
\Gamma $ is the group law.

The sheaf of $\tau$-functions of a point $U\in \gr^0(V)$ (not necessarily a
geometric point),
$\tilde {\L_\tau}(U)$, is the invertible sheaf over $\Gamma\times
\{U\}$ defined by:
$$\tilde{\L_\tau}(U)={\mathcal P}\vert_{\Gamma\times\{U\}}$$

The restriction homomorphism induces the following homomorphism
between global sections:
\beq
H^0\left(\Gamma\times \gr^0(V), \mu^*\det_V^*\right)\to
H^0\left(\Gamma\times \{U\},\tilde {\L_\tau}(U)  \right)
\label{eq:restriction}
\end{equation}

\begin{defn}\label{defn:tau0}
The $\tau$-function of the point $U\in\gr^0(V)$ over $\Gamma$ is defined
as the image $\tilde {\tau }_U$ of the section $\mu^*\Omega_+$ by the
homomorphism {\ref{eq:restriction}} ($\Omega_+$ being the global section
defined in {\ref{defn:global-section}}).
\end{defn}

Obviously $\tilde {\tau }_U$ is not a function over $\Gamma\times
\{U\}$ since the invertible sheaf $\tilde {\L_\tau}(U)$ is not
trivial, but this definition is essentially the $\tau$-function defined by
M. and Y. Sato and Segal-Wilson (\cite{SS,SW}). Their
$\tau$-function is obtained by restricting the invertible sheaf
$\tilde {\L_\tau}(U)$ to the subgroup $\Gamma_-\subset \Gamma$.

From \cite{AMP} we know that the invertible sheaf over $\Gamma_-$:
$${\L_\tau}(U)= \tilde {\L_\tau}(U)\vert_{\Gamma_-\times \{U\}}$$
is trivial, and that:
$$\sigma_0(g)=g\cdot\delta_U$$
 (where $g\in
\Gamma_-)$ and $\delta_U$ be a non-zero element in the fibre of
${\L_\tau}(U)$ over the point $(1,U)$ of $\Gamma\times\{U\}$) is a
global section of ${\L_\tau}(U)$ without zeroes; and, therefore, it
gives a trivialization.

With respect to this, the global section of  ${\L_\tau}(U)$ defined
by $\tilde{\tau}_U$ is identified with the function
$\tau_U\in \o(\Gamma_-)=k\{\{x_1,\dots\}\}$ given by Segal-Wilson
\cite{SW}:
$${\tau}_U(g)=\frac{\tilde{\tau}_U}{\sigma_0}=
\frac{\mu^*\Omega_+}{\sigma_0}=\frac{\Omega_+(gU)}{g\cdot\delta_U}$$

Observe that the $\tau$-function ${\tau}_U$ is not a series
of infinite variables but rather  an element of the ring $k\{\{x_1,
\dots\}\}$.

\begin{rem}
In the literature (\cite{SS,AD,KNTY}) one
also finds  another definition of the $\tau$-function:  for a
geometric point $\tilde U\in\det_V$ in the fibre of $U\in\gr^0(V)$, the
$\bar\tau$-function of
$\tilde U$ is defined as the element
$\bar \tau(U)\in H^0\left (\detd_V\right)^*\otimes k(U)$ ($k(U)$
being the residual field of $U$).

The deep relationship between both definitions emerges through the so
called boso\-ni\-zation isomorphism. To  introduce this isomorphism
certain preliminaries are necessary.

Since the subgroup $\Gamma_+$ of $\Gamma$ acts
freely on $\gr^0(V)$, the orbits of the rational points of $\gr^0(V)$
under the action of $\Gamma_+$ are isomorphic to $\Gamma_+$ (as schemes).
Let $X$ be the orbit of $V_-=z^{-1}\cdot k[z^{-1}]\subset V$ under
$\Gamma_+$. The restrictions of $\det_V$ and $\detd_V$ to $X$ are
trivial invertible sheaves. Bearing in mind that the points of
$X$ are $k$-vector subspaces of $V$ whose intersection with $V_+$ is
zero, one has that the section $\Omega_+$ of $\det_V^*$ defines a
canonical trivialization of $\detd_V$ over $X$.

Now, the {\bf bosonization isomorphism} is the canonical isomorphism:
$$\Omega(S)\longrightarrow \o(\Gamma_+)$$
($\Omega(S)$ is defined in {\ref{defn:global-section}}) induced by the
restriction homomorphism:
$$B: H^0(\gr^0(V),\detd_V)\rightarrow  H^0(X,\detd_V\vert_X)$$
and the isomorphism:
$$ H^0(X,\detd_V\vert_X)\iso\o(\Gamma_+)$$
 (associated to the trivialization defined by $\Omega_+$).

Finally, note that there exists an isomorphism of $k$-vector
spaces:
$$\o(\Gamma_+)^*=k[x_1,\dots]^*\to
\o(\Gamma_-)=k\{\{x_1,x_2,\dots\}\}$$
 identifying the Schur polynomial $F_S$ with the linear form:
$$F_{S'}\mapsto F_S(F_{S'})=\delta_{S,S'}$$
 (for details see the first chapter of
\cite{Mc}). Now, the composition of the homomorphism
$B^*$ (the dual of $B$) and the above isomorphism gives:
$$\tilde B^*:\o(\Gamma_+)^*=k\{\{x_1,x_2,\dots\}\}
\longrightarrow H^0\left (\detd_V\right)^*$$
The relationship between $\tau_U$ and $\bar \tau (\tilde U)$ is:
$\tilde B^*(\tau_U)=\bar \tau (\tilde U)$ (up to a non-zero constant).

The transformations of vertex operators under this isomorphism can now
be explicitely computed when the characteristic is zero (\cite{DJKM}).
In our approach, vertex operators are to be understood as the formal
Abel morphisms which will be studied below.
\end{rem}

Once we have defined the $\tau$-function algebraically, we can define
the Baker-Akhiezer functions using formula~5.14. of \cite{SW}; this
procedure is used by several authors. However, we prefer to continue with
the analogy with the classical theory of curves and Jacobians and define
the Baker-Akhiezer functions as a formal analogue of the universal invertible
sheaf of the Jacobian.

Let us consider the composition of morphism:
$$\tilde \beta:\hat C\times \Gamma\times \gr^0(V)
\overset {\phi\times Id}\longrightarrow
\Gamma \times\Gamma\times\gr^0(V)\overset {m\times Id}
\longrightarrow\Gamma\times \gr^0(V)$$

$\phi:\hat C=\sf k[[z]] \to \Gamma$ being the Abel morphism (taking
values in $\Gamma_-\subset \Gamma$) and
$m:\Gamma\times \Gamma\to\Gamma$ the group law.

\begin{defn}
The sheaf of Baker-Akhiezer functions is
the invertible sheaf over $\hat C\times \Gamma\times \gr^0(V)$ defined by:
$$\tilde\L_{\tiny BA}=(\phi\times Id)^*(m\times Id)^*{\mathcal P}$$

Let us define the sheaf of Baker-Akhiezer functions at  a point $U\in
\gr^0(V)$ as the invertible sheaf:
$$\tilde\L_{\tiny BA}(U)=\tilde\L_{\tiny BA}\vert_{\hat C\times
\Gamma\times \{U\}}=\tilde {\beta_U}^*\tilde {\L_\tau}(U)$$
where $\tilde {\beta_U}^*$ is the following homomorphism
between global sections:
$$ H^0(\Gamma\times \{U\},\tilde {\L_\tau}(U))
\overset{\tilde{\beta_U}^*}\longrightarrow
 H^0(\hat C\times\Gamma\times \{U\},\tilde\L_{\tiny BA}(U))$$
\end{defn}

By the definitions,
$\tilde\L_{\tiny BA}(U)\vert_{\hat C\times\Gamma_-\times
\{U\}}$  is a trivial invertible sheaf over $\hat C \times
\Gamma_-$.

\begin{defn}
The Baker-Akhiezer function of a point $U\in \gr^0(V)$ is
 $\psi_U=v^{-1}\cdot \beta^*_U(\tau_U)$ where $\beta^*$ is the
restriction homomorphism:
$$H^0\left( \Gamma\times \{U\},\tilde
{\L_\tau}(U)\right)\overset{\beta^*}\longrightarrow
 H^0\left(\hat C\times
\Gamma\times \{U\},\tilde\L_{\tiny BA}(U)\right)$$
induced by $\tilde{\beta_U}^*$.
\end{defn}

Note that from the definition one has the following
expression for the Baker-Akhiezer function:
\beq
\psi_U(z,g)=v(g,z)^{-1}\cdot \frac{\tau_U\left(g\cdot
\phi_1(z)\right)}{\tau_U(g)}
\label{eq:BA}\end{equation}
and that the Baker-Ahkiezer function of $V_-=z^{-1}\,k[z^{-1}]$ is the
universal invertible element $v^{-1}$ (see remark {\ref{rem:uni-ele}}).

When the characteristic of the base field $k$ is zero, we can
identify $\Gamma_-$ with the additive group scheme
 ${\hat\A}_{\infty}$ through the exponential. Therefore, the latter
expression
 is the classical expression for the Baker-Akhiezer functions
(\cite{SW}~5.16):
$$\psi_U(z,t)=\exp(-\sum t_i\,z^{-i})\cdot
\left(\frac{\tau_U(t+[z])}{\tau_U(t)}\right)$$
where $[z]=(z,\frac{1}{2}z^2,\frac{1}{3} z^3,\dots)$ and
$t=(t_1,t_2,\dots)$ and, through the
exponential map, $v(t,z)=\exp(\sum t_i\,z^{-i})$.

For the general case, we obtain explicit expressions for $\psi_U$ as a
function over
$\hat C\times {\hat\A}_{\infty}$ but considering in
${\hat\A}_{\infty}$ the group law, $*$, induced by the exponential
{\ref{defn:exp-gamma-minus}}:
$$\psi_U(z,t)=v(t,z)^{-1}\cdot \frac{\tau_U\left(t*
\phi(z)\right)}{\tau_U(t) }$$

\begin{rem}
Note that our definitions of
$\tau$-functions and Baker-Akhiezer functions are valid over arbitrary
base schemes. One then has the notion of $\tau$-functions and
Baker-Akhiezer  functions for families of elements of $\gr^0(V)$ and, if
we consider the Grassmannian of $\Z((z))$ one then has $\tau$-functions
and Baker-Akhiezer functions of the rational points of
$\gr\left(\Z((z))\right)$ and the geometric properties studied in this
paper have a translation into arithmetic properties of the elements of
$\gr\left(\Z((z))\right)$. The results stated by Anderson in \cite{An}
are a particular case of a much more general setting valid not only for
$p$-adic fields but also for arbitrary global numbers field.
\end{rem}

The fundamental property of the $\tau$-function is the
analogue of the Addition Formulae.

Let $\phi_N$ be the Abel morphism of degree $N$ ($N>0$); that is, the
morphism $\hat C^N\to\Gamma_-$ given by
$\prod_{i=1}^N(1-\frac{x_i}z)^{-1}$, where $\hat C^N=\sf(A)$
($A=k[[x_1,\dots,x_N]]$). Let $\Gamma_-=\sf k\{\{t_1,\dots\}\}$ and
observe that, in this setting, $t_i$ should be
interpreted as the coefficient of $z^{-i}$ in
$\prod(1-\frac{x_j}z)^{-1}$.

For simplicity's sake, $U_A$ will denote the
point $U\hat\otimes A\in\gr^0(V)(A)$ for a rational point
$U\in\grv$, and $V$ (resp. $V_+$) will denote $\w V_A$
(resp. $(\w{V_+})_A$).

Now, for our $\tau$-function we shall give formal analogues of the
addition formula of \cite{SS} and of corollary~2.19 of
\cite{F} for theta functions and of Lemma~4.2 of \cite{Kon} for
$\tau$-funtions. Let us begin with an explicit computation for $\tau$.

\begin{lem}\label{lemma:add-for}
Let $U$ be a rational point of $\gr^0(V)$.
Assume that $V/V_++z^N\cdot U=0$ and $V_+\cap z^N\cdot U$ is
$N$-dimensional, and let
$\{f_1,\dots,f_N\}$ be a basis of the latter. One then has that:
$$\phi_N^*\tau_U = \frac{1}{\prod_{i< j}(x_i-x_j)}\cdot
det \pmatrix
f_1(x_1) & \dots & f_1(x_N) \\
\vdots & &\vdots \\
f_N(x_1) & \dots & f_N(x_N)
\endpmatrix$$
as functions on $\hat C^N$ (up to elements of $k^*$).
\end{lem}

\begin{pf}
By the very definition of the $\tau$-function and by the properties of
$Det$, it follows that $\phi_N^*\tau_U$ equals the determinant of the
inverse image of the complex:
$$\L\to\left({ V}/{ V_+}\right)_{\gr^0(V)}$$
by the morphism $\hat C^N\to\Gamma_-\to\gr^0(V)$, which is precisely:
$${\mathcal
C}^\punto\equiv g\cdot U_A\to{A((z))}/{A[[z]]}=V/V_+$$

Let us define the following homomorphism of $A$-modules:
$$\aligned \alpha_N:A[[z]] & \to A^N \\
f(z)&\mapsto\left(f(x_1),\dots,f(x_N)\right)\endaligned$$
whose kernel is the ideal generated by $\prod_{i=1}^N (z-x_i)$. One
thus has the following exact sequence of complexes of $A$-modules:
{\small $$\minCDarrowwidth17pt\CD
0 @>>> \prod(z-x_i)\cdot { V_+} @>>> { V_+} @>>>
{{ V_+}}/{\prod(z-x_i)\cdot { V_+}} @>>> 0 \\
@. @V{\beta}VV @V{(\pi,\alpha_N)}VV @V{\bar\alpha_N}VV \\
0 @>>> V/{z^N\cdot U_A} @>>>
\left(V/{z^N\cdot U_A}\right)\oplus A^N  @>>> A^N  @>>> 0
\endCD$$ }

The complex on the middle (right hand side resp.) will be denoted by
${\mathcal C}^\punto_1$ (${\mathcal C}^\punto_2$ resp.). Further, note
that the complex on the left hand side is quasi-isomorphic to ${\mathcal
C}^\punto$.

\noindent Observe that:
\begin{itemize}
\item $\det({\mathcal C}^\punto_2)$ is isomorphic to the ideal of the
diagonals (as an $A$-module), and the section
$det(\bar\alpha_N)$ is exactly $\prod_{i<j}(x_i-x_j)$,
\item $det(\beta)=det((\pi,\alpha_N))\cdot det(\bar\alpha_N)^{-1}$.
\end{itemize}

Moreover, we also have another exact sequence:
{\small $$\minCDarrowwidth17pt\CD
0 @>>> { V_+}\cap z^N\cdot U_A  @>>>
{ V_+} @>>> {{ V_+}}/{{ V_+}\cap z^N\cdot U_A} @>>>0 \\
@. @V{\alpha^U_N}VV @V{(\alpha_N,\pi)}VV @V{\simeq}VV \\
0 @>>> A^N  @>>> A^N  \oplus{ V}/{z^N\cdot U_A} @>>>
{ V}/{z^N\cdot U_A} @>>> 0
\endCD$$ }
Let ${\mathcal C}^\punto_3$ denote the complex on the right hand side.
The hypothesis implies that the non-trivial homomorphism of the complex
${\mathcal C}^\punto_3$ is an isomorphism. They also imply that
$\det({\mathcal C}^\punto_3)\in k^*$.

From both these  exact sequences one has the following relation:
$$det(\beta)=\frac{det(\alpha^U_N)}{\prod_{i<j}(x_i-x_j)}\qquad \text{(up to
an element of $k^*$)}$$

In order to compute $det(\alpha^U_N)$, let us choose $\{f_1,\dots,f_N\}$
a basis of ${ V_+}\cap z^N\cdot U_A$ (as an $A$-module). The
matrix associated to $\alpha^U_N$ is now:
$$\pmatrix
f_1(x_1) & \dots & f_1(x_N) \\
\vdots & &\vdots \\
f_N(x_1) & \dots & f_N(x_N)
\endpmatrix$$
and the formula follows.
\end{pf}

This proof implies that (in the same conditions) $\phi_N^*\tau_{hU}=\prod_i
h(x_i)\cdot\phi_N^*\tau_U$, where $h\in\Gamma_+$, and also the following
theorem.

\begin{thm}[Addition Formulae \cite{SS}]
Let $U\in\gr^0(V)$ be a rational point; then, for all $n\leq N$ (natural
numbers) the functions:
$$\left\{ \prod_{j<k}(x_{i_k}-x_{i_j}) \cdot \phi_{i_1\dots i_n}^*\tau_U
\right\}_{0<i_1<\dots<i_n\leq N}$$
satisfy the Pl\"ucker equations. Here, $\phi_{i_1\dots i_n}$ denotes the
morphism
$\hat C^N\to\Gamma_-$ given by $\prod_j(1-\frac{x_{i_j}}z)^{-1}$.
\end{thm}

We finish this section with a characterization of the Baker-Akhiezer
function. The importance of this result will be clear in the next
section.
Further, this characterization will show the close relationship between
the Baker-Akhiezer of a point $U\in\gr^0(V)$ and a basis of $U$ as
$k$-vector space. (Compare with Proposition 5.1 of \cite{SW} and
Proposition 4.8 of \cite{KNTY}).

\begin{thm}\label{4:thm:BA}
Let $U\in\gr^0(V)$ be a rational point.
Then:
$$\psi_U(z,t)=
z\cdot\sum_{i\geq 1}
\psi_U^{(i)}(z) p_i(t)$$
where $\psi_U^{(i)}(z)\in U$ and $p_i(t)\in k\{\{t_1,\dots\}\}$ is
independent of $U$.

Furthermore, let
$U$ be in $F_{A_S}$ for a sequence $S$ and
$g$ be $\prod_{j=1}^n(1-\frac{x_j}z)^{-1}\in\Gamma_-$. Define $s_i$ by
$\Z-S=\{s_1,s_2,\dots\}$ (as in {\ref{subsect:comp-inf-grass}}). Then
$\psi_U^{(i)}(z)$ has a pole in
$z=0$ of order $s_i$, and $p_i(x)$ is a homogeneous polynomial in the
$x$'s of degree $i-1$.
\end{thm}

\begin{pf}
Since $\Gamma_-=\sf k\{\{t_1,\dots\}\}$ is the direct limit of the
symmetric products of $\hat C$ (recall the relation between the $x$'s
and $t$'s), let us compute $\psi_U\vert_{\hat C^N}$ using
{\ref{lemma:add-for}} for $N>>0$. Choose $N$ such that
$z^{-N-i}k[[z]]\cap U$ is $N+i$-dimensional for $i=0,1$, and let
$\{f_1(z),\dots,f_N(z)\}$ be a basis of $z^{-N}k[[z]]\cap U$, and
$\{f_1(z),\dots,f_{N+1}(z)\}$ of
$z^{-N-1}k[[z]]\cap U$. With no loss of generality, we can assume that
$f_i$ has a pole at $z=0$ of order $s_i$.

Recall that for $g=\prod_{j=1}^n(1-\frac{x_j}z)^{-1}$ the Baker-Akhiezer
function is:
\beq
\psi_U(z,g)=g(z)^{-1}\cdot
\frac{\tau_U(g\cdot\phi(z))}{\tau_U(g)}
\label{eq:eq1}
\end{equation}

Let $\phi_{1,N}$ be the morphism $\hat C\times\hat C^N\to\Gamma_-$ given
by $g\cdot(1-\frac{\bar z}z)^{-1}$, and $\phi_N$ the morphism $\hat
C^N\to\Gamma_-$ given by $g$. Here, the ring of $\hat C^N$ is
$k[[x_1,\dots,x_N]]$ and the ring of $\hat C\times\hat C^N$ is
$k[[\bar z]]\hat\otimes k[[x_1,\dots,x_N]]$.

By lemma {\ref{lemma:add-for}} one sees that:
$$\frac{\phi^*_{1,N}\tau_U}{\phi^*_N\tau_U}
=(-1)^N\cdot\prod_{i=1}^N(\bar z-x_i)^{-1}\cdot
\left(\bar f_{N+1}(\bar z)\prod_{i}x_i+
\sum_{i=1}^N\bar f_i(\bar z)p_i(x_1,\dots,x_N)\right)$$
where $\bar f_i(z)=z^{N+1}f_i(z)$ and $p_i$ is a symmetric polynomial.
Now,  from the above expression for the Baker-Akheizer function, one has:
$$\psi_U(z,t)\vert_{\hat C^N}=
z\cdot\sum_{i=1}^{N+1}f_i(z)p_i(x_1,\dots,x_N)$$
and therefore $\psi_U(z,t)= z\cdot\sum_{i>0}f_i(z)p_i(t)$, where
$f_i(z)\in U$ and $p_i(t)\in k\{\{t_1,\dots\}\}$. Further, since
$\{\psi_U(z,t)\vert_{\hat C^N}\}_{N\in{\mathbb N}}$ is an element of an
inverse limit and $p_{N+1}(x)\vert_{\hat C^N}$ is known, one can compute
$p_i$ explicitly. From the choice of $f_i$ and the properties of
$p_i$ one concludes.
\end{pf}

Now, all the above definitions and results on tau and BA functions can be
generalized to the case of $U\in\gr^n(V)$, a point in an arbitrary
connected component.

\begin{defn}[$\tau$ and BA functions for arbitrary points]
Let $U$ be a point of $\gr^n(V)$ ($n\in\Z$). Define its $\tau$-function
by Definition \ref{defn:tau0} but replacing $\Omega_+$ by $\Omega_+^n$
(see Remark \ref{rem:indexn}). Define its Baker-Akhiezer function by
formula \ref{eq:BA}.
\end{defn}

Observe that from the definition of $\tau$-function, we have that:
$$\tau_U(g)\,=\,\frac{\Omega_+^n(gU)}{g\cdot \delta_U}\,=\,
\frac{\Omega_+(g\cdot z^{-n}U)}{g\cdot \delta_U}$$
Furthermore, Lemma \ref{lemma:add-for} can be generalized in the
following way:
let $U$ be a rational point of $\gr^n(V)$.
Assume that $V/V_++z^{N-n}\cdot U=0$ and $V_+\cap z^{N-n}\cdot U$ is
$N$-dimensional, and let
$\{f_1,\dots,f_N\}$ be a basis of the latter. One then has that:
$$\phi_N^*\tau_U = \frac{1}{\prod_{i< j}(x_i-x_j)}\cdot
det \pmatrix
f_1(x_1) & \dots & f_1(x_N) \\
\vdots & &\vdots \\
f_N(x_1) & \dots & f_N(x_N)
\endpmatrix$$
as functions on $\hat C^N$ (up to elements of $k^*$).

Similarly, the very definition of Baker-Ahkiezer function implies that:
$$\psi_U(z,t)\,=\,\psi_{z^{-n}U}(z,t)$$
for a point $U\in\gr^n(V)$. Moreover, the proof of Theorem \ref{4:thm:BA}
shows that in this case:

\beq
\psi_U(z,t)\,=\,z^{1-n}\cdot\sum_{i\geq 1}
\psi_U^{(i)}(z) p_i(t)
\label{eq:BA=sum}\end{equation}
where $\psi_U^{(i)}(z)\in U$ and $p_i(t)\in k\{\{t_1,\dots\}\}$ is
independent of $U$.

\section{Bilinear Identity and KP Hierarchy}

This section has two well defined parts. In the first part the famous
Residue Bilinear Identity is deduced, while in the second one the
equivalence with the KP hierarchy is shown, which was already known
(see \cite{DJKM}). Nevertheless, the importance lies not in the result
but in the proofs, since the methods used will also allow us to prove the
fundamental theorems of our paper.

The essential ingredient comes from the relationship between
$\gr(k((z)))$ and $\gr(k((z))^c)$ outlined in {\ref{subsect:grass-dual}}.
However, in our case there exists a metric on $V$; namely, that induced
by the residue pairing:
$$(f,g) \,=\, \res_{z=0} (f\cdot g) dz $$

Let us denote by $\bar z^i$ the element of $V^c$ such that $\bar
z^i(z^j)=\delta_{ij}$. One can formally write $V^c=k((\bar z^{-1}))$ (as
$k$-vector spaces), and it may be seen that  the homomorphism induced by
the residue:
$$\aligned V&\to V^c \\
z^i &\mapsto \bar z^{-i-1} \endaligned$$
is in fact an isomorphism, and sends $V_+$ to $V_+^\diamond$. It
therefore induces  an isomorphism
$\gr(V^c,\B^\diamond)\iso\gr(V,\B)$.

The composition of the latter isomorphism and the isomorphism $I$
constructed in {\ref{subsect:grass-dual}} gives a non-trivial automorphism
of the
Grassmannian:
$$\aligned R:\grv &\longrightarrow \grv \\
L &\longmapsto L^{\perp} \endaligned$$

Trivial calculation shows that $R^*\det_V\simeq \det_V$, and that the
index of a point $L\in\grv$ is exactly the opposite of the index of
$R(L)=L^\perp$. Therefore, this induces an involution:
$$R^*: H^0(\gr^n(V),\detd_{n})\to H^0(\gr^{-n}(V),\detd_{-n})$$
It is now straightforward to see that $R^*\Omega_+^{-n}=\Omega_+^{n}$.

The last remarkable fact about $R$ is that for a given point $U\in\grv$ the
morphism:
$$\mu_U:\Gamma\to \grv$$
($\mu_U$ being that induced by the action of $\Gamma$ on $\grv$) is
equivariant with respect to ``passing to the inverse'' in $\Gamma$ and
$R$ in $\grv$. In other words:
$$(g\cdot U)^\perp= g^{-1}\cdot U^\perp$$

From the latter two observations, one has trivially $\Omega_+^n(g\cdot
U^\perp)=\Omega_+^n((g^{-1}\cdot U)^\perp)=\Omega_+^{-n}(g^{-1}\cdot U)$,
and hence:
$$\tau_{U^\perp}(g)=\tau_U(g^{-1})$$
This motivates us to give the following:

\begin{defn}
The adjoint Baker-Akheizer function of a point $U\in\grv$ is:
$$\psi^*_U(z,g)=\psi_{U^\perp}(z,g^{-1})$$
(Recall that: $\psi_{U^\perp}(z,g)=
g^{-1}\frac{\tau_U(g^{-1}*\phi(z)^{-1})}{\tau_U(g^{-1})}$).
\end{defn}

Note that in characteristic zero, the definition gives:
$$\psi^*_U(z,t)=\exp(\sum_{i>0}t_i\,z^{-i})\frac{\tau_U(t-[z])}{\tau_U(t)}$$

By the very definition of the Baker-Akhiezer function and formula
\ref{eq:BA=sum}, it follows that:
$$\left(\psi_{U^\circ}(\bar z^{-1},t')\right)
\left(\psi_U(z,t)\right) \,=\,0$$
or, in other words:
$$\res_{z=0} \psi_U(z,t)\psi_{U^\perp}(z,t')\frac{dz}{z^2}\,=\,0$$

Using the adjoint Baker-Akhiezer function of $U$ as $\psi_{U^\perp}$, then the
latter equation can be rewritten as:
\beq
\res_{z=0}
\psi_U(z,t)\psi^*_{U}(z,t')\frac{dz}{z^2}\,=\,0
\label{eq:residue}
\end{equation}

\begin{thm}\label{thm:res-equa}
Let $U$ and $U'$ be two rational points of the Grassmannian and let us
assume they have the same index. Then, it holds that:
$$\res_{z=0}
\psi_U(z,t)\psi^*_{U'}(z,t')\frac{dz}{z^2}\,=\,0\quad\iff\quad\ U=U'$$
\end{thm}

\begin{pf}
The converse is already shown. For the direct one, observe that the
identity and formula \ref{eq:BA=sum} imply that: ${U'}^\perp\subseteq
U^\perp$. Recalling that both have the same index and that the metric is
non-degenerate, the result follows.
\end{pf}

When the base field is the field of complex numbers, $\mathbb C$, it is
well known that the above results (valid over arbitrary fields) turn out
to be equivalent to the KP hierarchy (see \cite{DJKM,F2}).

The goal now is to show how the KP hierarchy is obtained
from the Residue Bilinear Identity (when $k={\mathbb C}$). In this
procedure we shall use Schur polynomials and their properties in a very
fundamental way. However, the advantage of this is that we can deal in a
similar way with the residue condition that characterizes the moduli
space of pointed curves in the Grassmannian, and this will allow us to
compute its equations explicitly.

Let us begin with some preliminaries on symmetric polynomials following
the first chapter of \cite{Mc}. Let $k$ be a ring of characteristic zero,
and let
$\mathcal S$ denote the subring of $k\{\{x_1,\dots\}\}$ consisting of symmetric
polynomials. Given a decreasing sequence $\lambda$  of natural numbers
$\lambda_1\geq\dots\geq\lambda_n$ define the associated Schur polynomial as:
$$\chi_\lambda(x)=\frac{det(x_j^{\lambda_i+n-i})}{det(x_j^{n-i})}$$

It is also known that there exists a non-degenerated metric  in $\mathcal
S$ for which the Schur polynomials are an orthonormal basis. Denote by
$\chi^*_\lambda$ the linear form ${\mathcal S}\to k$ defined by:
$$\chi^*_\lambda(\chi_\mu(x))=\delta_{\lambda,\mu}$$

Therefore, the identity morphism $Id:{\mathcal S}\to {\mathcal S}$ can be
expressed as $\sum_\lambda\chi_\lambda(x) \chi^*_\lambda$.
That is, for each element $p(x)\in{\mathcal S}$ one has that:
$$p(x)=\sum_\lambda \chi_\lambda(x)\chi^*_\lambda(p)$$

\begin{rem}\label{rem:Taylor}
\begin{enumerate}
\item Note that if $k$ is a field, then one can express $\chi^*_\lambda$ in
terms
of differential operators. Explicitly, this is:
$$\chi^*_\lambda=\chi_\lambda(\tilde\partial_x)\vert_{x=0}$$
where $\tilde\partial_x=
(\frac{\partial}{\partial x_1},\frac12\frac{\partial}{\partial
x_2},\dots)$. The operator $\sum_\lambda
\chi_\lambda(x)\chi_\lambda(\tilde\partial_x)\vert_{x=0}$ will be called
Taylor expansion operator. Analogously, for more than
one set of variables, for instance $x=(x_1,x_2,\dots)$ and
$y=(y_1,y_2,\dots)$, the identity can be written as:
$$Id=\sum_{\lambda,\mu}\chi_\lambda(x)\chi_\mu(y)
\chi_\lambda(\tilde\partial_x)\chi_\mu(\tilde\partial_y)\vert_{x=y=0}$$
\item Assume we have two set of variables, as above. Observe
that (applying the Taylor expansion operator to $f$):
$$\left(\sum_\lambda\chi_\lambda(y)\chi^*_\lambda(\tilde\partial_x)\right)f(x)
\vert_{x=0}=f(y)$$
That is, this operator replaces $x_i$ by $y_i$.
\item Similarly, by replacing $x_i$ by $\frac{z^i}i$
and denoting by $p_j$ the polynomials defined by $\exp(\sum_{i>0} x_i
z^i)=\sum_{j\geq 0} p_j(x)z^j$, one has:
$$\text{coefficient of $z^j$ in $f(z)$ }=\,
p_j(\tilde\partial_x)\vert_{x=0}f(x)$$
(note that $p_j(\tilde\partial_x)\vert_{x=0}(p_i(x))=\delta_{ij}$).
\item Now, observe that:
$$P(\partial_y)\vert_{y=0}f(x+y)=P(\partial_x)f(x)$$
where $P$ is a polynomial, $f$ a function, and $x+y$ denotes
$(x_1+y_1,x_2+y_2,\dots)$. It therefore follows that:
$$\text{coefficient of $z^j$ in }f(x+[z])\,=\,
p_j(\tilde\partial_y)\vert_{y=0}f(x+y) \,=\,
p_j(\tilde\partial_x)f(x)$$
(where $[z]$ is $(z,\frac{z^2}2,\dots)$), and hence:
$$f(x+[z])=\sum_j z^jp_j(\tilde\partial_x)f(x)$$
That is, the operator $\sum_j z^jp_j(\tilde\partial_x)$ (also written as
$\exp(\sum z^i\tilde\partial_{x_i})$) replaces $x_i$ by $x_i+\frac{z^i}i$.
Finally, this operator relates the $\tau$-function and the Baker-Akhiezer
function by:
$$\psi_U(z,t)=\exp(-\sum t_iz^{-i})\frac
{\exp\left(\sum z^i\tilde\partial_t\right)\tau_U(t)}{\tau_U(t)}$$
\end{enumerate}
\end{rem}

We can now state the main result of this subsection in the form given in
\cite{F2}.

\begin{thm}[KP Equations]\label{thm:KP-eq}
The condition:
$$\operatorname{Res}_{z=0}\psi_U(z,t)\cdot\psi^*_U(z,t')\frac{dz}{z^2}=0$$
for a rational point $U\in\grv$ is equivalent to the infinite set of
equations (indexed by a pair of Young diagrams $\lambda_1,\lambda_2$):
$$\left(\sum
p_{\beta_1}(\tilde\partial_t)D_{\lambda_1,\alpha_1}(-\tilde\partial_t)
\cdot
p_{\beta_2}(-\tilde\partial_{t'})D_{\lambda_2,\alpha_2}(\tilde\partial_{t'})
\right)\vert_{t=t'=0}\tau_U(t)\cdot \tau_U(t')\,=\,0$$
where the sum is taken over the 4-tuples
$\{\alpha_1,\beta_1,\alpha_2,\beta_2\}$ of integers such that
$-\alpha_1+\beta_1-\alpha_2+\beta_2=1$, and
$D_{\lambda,\alpha}=\sum_\mu \chi_\mu$ where $\mu$ is a Young diagram
such that $\lambda-\mu$ is a horizontal $\alpha$-strip.
\end{thm}

\begin{pf}
First, observe that $\Gamma_-$ (and equivalently $\Gamma_+$) acts on
$\grv$ by homotheties preserving the determinant sheaf. Then, by
straightforward computation one shows that
$\tau_U(t)=\sum_{\lambda}\Omega_\lambda(U)
\chi_\lambda(t)$ (the sum taken over the set of Young diagrams).

Now recall two basic facts about functions  $f(x)\in
k\{\{x_1,x_2,\dots\}\}$: the coefficient of $z^\beta$ in $f([z])$
is precisely $p_\beta(\tilde\partial_y)f(y)\vert_{y=0}$ ($[z]$ being
$(z,\frac12z^2,\dots)$); and
$\chi_\lambda(\tilde\partial_y)\vert_{y=0}f(x+y)=
\chi_\lambda(\tilde\partial_x)f(x)$.
From them it follows that:
$$\text{coefficient of $z^\beta$ in }f(x+[z])\,=\,
p_\beta(\tilde\partial_x)f(x)$$

Let us now begin, properly speaking, with the proof of the theorem. The
residue condition is trivially equivalent to:
$$\text{coefficient of $z$ in }
\exp(-\sum t_iz^{-i})\tau_U(t+z)\exp(\sum t'_iz^{-i})\tau_U(t'-z)\,=\,0$$
This coefficient is given by the sum:
$$\sum p_{\alpha_1}(-t)p_{\beta_1}(\tilde\partial_t)\tau_U(t)\cdot
p_{\alpha_2}(t')p_{\beta_2}(-\tilde\partial_{t'})\tau_U(t')$$
over the 4-tuples $\{\alpha_1,\beta_1,\alpha_2,\beta_2\}$ of integers such
that $-\alpha_1+\beta_1-\alpha_2+\beta_2=1$.
But now a necessary and
sufficient condition for an element $f(t,t')\in k\{\{t,t'\}\}$ to be zero is
that $\chi_{\lambda_1}(\tilde\partial_t)\chi_{\lambda_2}(\tilde\partial_{t'})
\vert_{t=t'=0}$ applied to $f$ must be zero for all pairs of
Young diagrams $\lambda_1,\lambda_2$.
Simple computation together with formula {\bf I}.5.16 of \cite{Mc} gives
the following infinite set of quadratic differential equations
on $\tau_U$:
$$\left(\sum
p_{\beta_1}(\tilde\partial_t)D_{\lambda_1,\alpha_1}(-\tilde\partial_t)
\cdot
p_{\beta_2}(-\tilde\partial_{t'})D_{\lambda_2,\alpha_2}(\tilde\partial_{t'})
\right)\vert_{t=t'=0}\tau_U(t)\cdot \tau_U(t')\,=\,0$$
where $D_{\lambda,\alpha}=\sum_\mu\chi^*_\mu$
(the sum taken over the set of Young diagrams $\mu$ such that
$\lambda-\mu$ is a horizontal $\alpha$-strip).
\end{pf}

\begin{rem}
Recalling Remark {\ref{rem:grass=plucker}} one obtains four different
ways to characterize the set of rational points of the infinite
Grassmannian $\gr({\mathbb C}((z)))$ into the infinite dimensional
projective space
$\P\Omega(S)^*$; namely, the Pl\"ucker equations (\cite{Pl2}), the
Bilinear Residue Identity (see Proposition 4.15 of \cite{KNTY}, the
KP hierarchy (\cite{DJKM}) and the Hirota's Bilinear equations
(\cite{DJKM}). In the last two characterizations the differential
operators introduced in Remark \ref{rem:Taylor} are needed.
\end{rem}

\section{Characterization and Equations of the moduli space of pointed
curves in the Grassmannian}

In this section, and since our goal is to compute equations for the moduli
of pointed curves, we study a slightly modified Krichever construction;
namely, the application:
$$\aligned\{(C,p,\phi)\} &\to \grv \\
(C,p,\phi) &\mapsto \phi\left(H^0(C-p,\o_C)\right)
\endaligned$$
Here and in the sequel $\grv$ will denote the infinite Grassmannian of
the data $(V=k((z)),\B,V_+=k[[z]])$ as in Example \ref{2:exam}.

Let us introduce some more notation. Given a $k$-scheme $S$ define:
$$\begin{aligned}
\o_S[[z]] \,&=\, \limpl{n}\o_S[z]/z^n \\
\o_S((z)) \,&=\, \limil{m}z^{-m}\o_S[[z]]
\end{aligned}$$
Given a flat curve $\pi:C\to S$ and a section $\sigma: S\to C$ defining
a Cartier divisor $D$ denote:
$$\widehat\o_{C,D} \,=\, \limpl{n}\o_C/\o_C(-n)$$
where $\o_C(-1)$ is the ideal sheaf of $D$. Observe that
$\widehat\o_{C,D}$ is supported along $\sigma(S)$ and  it is then a
sheaf of $\o_S$-algebras. We also define the following sheaf of
$\o_S$-algebras:
$$\widehat \Sigma_{C,D}\,=\,\limil{m}\widehat\o_{C,D}(m)$$

\begin{defn}
Let $S$ be a $k$-scheme. Define the functor $\tilde\M^g_{\infty}$
over the category of $k$-schemes by:
$$S\rightsquigarrow \tilde\M^g_{\infty}(S)=
\{\text{ families $(C,D,\phi)$ over $S$ }\}$$
where these families satisfy:
\begin{enumerate}
\item $\pi:C\to S$ is a proper flat morphism,
whose geometric fibres are integral curves of arithmetic genus $g$,
\item $\sigma:S\to C$ is a section of $\pi$, such that when
considered as a Cartier Divisor $D$ over $C$ it is smooth, of relative
degree 1, and flat over
$S$. (We understand that $D\subset C$ is smooth over $S$, iff for
every closed point $x\in D$ there exists an open neighborhood $U$ of
$x$ in $C$ such that the morphism $U\to S$ is smooth).
\item $\phi$ is an isomorphism of $\o_S$-algebras:
$$\widehat \Sigma_{C,D}\,\iso\, \o_S((z))$$
\end{enumerate}
\end{defn}

On the set $\tilde\M^g_{\infty}(S)$ one can define an
equivalence  relation, $\sim$: $(C,D,\phi)$ and $(C',D',\phi')$ are said
to be equivalent, if there exists an isomorphism $C\to C'$ (over $S$)
such that the first family goes to the second under the induced
morphisms.

\begin{defn}
The moduli functor of pointed curves of genus $g$, $\M^g_{\infty}$, is the
functor over the category of
$k$-schemes defined by the sheafication of the functor:
$$S\rightsquigarrow {\tilde\M^g_{\infty}(S)}/{\sim}$$
(the superindex $g$ will be left out to denote the union over all
$g\geq 0$).
\end{defn}

\begin{prop}
The sheaf $\limil{m}\pi_*\o_{C,D}(m)$ is an $S$-valued point of
$\gr^{1-g}(\widehat \Sigma_{C,D})$ for all $(C,D,\phi)\in\M^g_{\infty}$.
\end{prop}

\begin{pf}
Consider the following exact sequence:
$$0\to \o_{C,D}(-n)\to \o_{C,D}(m)\to \o_{C,D}(m)/\o_{C,D}(-n)\to 0$$
for $n,m\geq 0$. Take now $\pi_*$ and recall that:
$$\pi_*\left(\o_{C,D}(m)/\o_{C,D}(-n)\right)\,\iso\,
\widehat \o_{C,D}(m)/\widehat\o_{C,D}(-n)$$
since it is concentrated on $\sigma(S)$. One then has the following
long exact sequence:
$$0\to\pi_*\o_{C,D}(-n)\to\pi_* \o_{C,D}(m)\to
\widehat\o_{C,D}(m)/\widehat\o_{C,D}(-n)\to
R^1\pi_*\o_{C,D}(-n)$$
and hence:
$$0\to\pi_*\o_{C,D}(-n)\to\limil{m}\pi_* \o_{C,D}(m)\to
\widehat\Sigma_{C,D}/\widehat\o_{C,D}(-n)\to
R^1\pi_*\o_{C,D}(-n)$$
Recalling now from \cite{Al} a Grauert type theorem that holds in this
case: there exists $\{U_i\}$ a covering of $S$ and
$\{m_i\}$ integers such that $\pi_*\o_C(m_i)_{U_i}$ is locally free of
finite type, and $R^1\pi_*\o_C(m_i)_{U_i}=0$; the result follows.
\end{pf}

Observe now that $\phi$ induces an isomorphism of infnite Grassmannians:
$$\gr^{1-g}(\widehat \Sigma_{C,D})\,\iso\,\gr^{1-g}(V)$$
in a natural way since there exists an integer $r$ such that:
$$z^{n+r}\o_S[[z]]\,\subseteq\,\phi(\widehat\o_{C,D}(n))\,\subseteq
\,z^{n-r}\o_S[[z]]\qquad\forall n$$

Summing up:
$$\phi\left(\limil{m}\pi_*\o_{C,D}(m)\right)$$
is a $S$-valued point of $\gr^{1-g}(V)$ for all
$(C,D,\phi)\in\M^g_{\infty}$. In other words, we have defined a morphism
of functors:
$$\aligned
K:\M_{\infty}&\longrightarrow \grv\\
(C,D,\phi)&\longmapsto \phi\left(\limil{n}\pi_*\o_C(n)\right)
\endaligned $$
which will be called the {\bf Krichever morphism}.(Note that $K$
considered for the $\sp({\mathbb C})$-valued points  is the usual
Krichever map, see \cite{K,PS,SW}).

Now let us state the following characterization of the image of
$K$, which is well known in the complex case:

\begin{thm}\label{thm:krich-char}
A point $U\in\grv(S)$ lies in the image of the Krichever morphism, if and
only if  $\o_S\subset U$ and $U\cdot U\subseteq U$ (where $\cdot$ is the
product of
$V$).
\end{thm}

\begin{pf}
Assume we have such a point $U$ of $\grv(S)$. Since $\M_{\infty}$ and
$\grv$ are sheaves, we can assume that $S$ is affine,
$S=\sp(B)$. Now, $U\subset B((z))$ is a sub-$B$-algebra and has a natural
filtration $U_n$ by the degree of $z^{-1}$. Let $\mathcal U$ denote the
associated graded algebra. It is not difficult to prove that
$C=\proj_B({\mathcal U})$ is a curve over $B$.

Let $I$ be the ideal sheaf generated by the elements $a\in{\mathcal U}$
such that the homogeneous localization ${\mathcal U}_{(a)}^0$ is
isomorphic to $\mathcal U$. Since $I$ is locally principal, it defines a
section $\sigma:S\to C$.

Finally, from the inclusion $U\subset \o_S((z))$ one easily deduces an
isomorphism $\phi:\widehat \Sigma_{C,D}\iso\o_S((z))$.

An easy calculation shows that this construction and the Krichever
morphism are the inverse of each other.
\end{pf}

\begin{thm}
$\M_\infty$ is representable by a closed subscheme of $\grv$.
\end{thm}

\begin{pf}
By the preceding characterization of
$\M_\infty$, it will suffice to recall the following result of \cite{Al}:

Let $A$ be a $k$-algebra, and $V$ a sheaf of $A$-modules over the
category of $A$-algebras. Let $M$, $M'$ be two
quasi-coherent subsheaves of $V$, such that $V/M'$ is isomorphic to a
inverse limit of finite type free $A$-modules. There then exists an
ideal $I\subseteq A$ such that every morphism
$f:A\to B$ with the property $M_B\subseteq M'_B$ (as subsheaves
of $V_B$) factorizes through $A/I$. (The subindex $B$ denotes the
canonically induced sheaf of $B$-modules over the category of
$B$-algebras).

Assuming this result, and since $\hat V_S/L$ is isomorphic
to a inverse limit of finite type $A$-modules for $L\in\grv(S)$, one
deduces that the conditions of {\ref{thm:krich-char}} are closed.

Let us now prove the claim. By the hypothesis one has that
$(V/M')_B\iso V_B/M'_B$ for all morphisms $A\to B$. Let $f$ be a
morphism $A\to B$, then the canonical inclusion $j:(M+M')/M'\to
V/M'$ gives:
$$j_B^i:((M+M')/M')_B @>>> (V/M')_B\iso \limpl{i} L_i\to L_i$$
where $L_i$ are finite type free $B$-modules.
One has that $M_B\subseteq M'_B$ if and only if $j_B^i$ is identically
zero for all $i$.

Recall that given a sub-$A$-module $\bar M$ of a finite type
free $A$-module $L=A^n$, the inclusion morphism $\bar j:\bar M\to L$
assigns to each $\bar m\in\bar M$ a set of coordinates $(\bar
j_1(\bar m),\dots,\bar j_n(\bar m))$, and hence $\bar j_B$ is identically
zero if and only if $f:A\to B$ factorizes through the ideal generated by
$\{j_i(\bar m)\,\vert\, \bar m\in \bar M, i=1,\dots, n\}$. This concludes
the proof.
\end{pf}

Our aim, now, is to give explicit characterizations of the set of
rational points $U$ satisfying the conditions of theorem
{\ref{thm:krich-char}}.

\begin{thm}
Let $S$ be a Young diagram of virtual cardinal $n$, and let $U$ be a
rational point of $F_{A_S}\subset\gr^n(V)$. The following three
conditions are equivalent:
\begin{enumerate}
\item $k\subset U$ and $U\cdot U\subseteq U$,
\item $U\cdot U=U$,
\item $0\notin S$ and
$\res_{z=0}\psi_U(z,t)\psi_U(z,t')\psi^*_U(z,t'')z^{n-3}dz=0$.
\end{enumerate}
\end{thm}

\begin{pf}
$1\implies2$ is trivial.

For $2\implies1$, one has only to check that $k\subset
U$. But recall that the element $u$ of $U-\{0\}\subset k((z))-\{0\}$ of
highest order is unique (up to a non-zero scalar) and should therefore
satisfy $u=\lambda\cdot u\cdot u$ and hence
$u=\lambda^{-1}\in k-\{0\}$.

$1\implies3$  First, note that $k\subset U$ implies
$0\notin S$. It is now clear by {\ref{eq:residue}} that the first
condition implies the third.

$3\implies1$ If the residue
condition is verified, it then follows that
$U^\perp\subseteq (U\cdot U)^\perp$ and therefore $U\cdot U\subseteq
U$, as desired. Now, since $0\notin S$ an element $u\in U$ of highest
order must have non-negative order, and since $u\cdot u\in
U$, one concludes that
$u=\lambda\in k-\{0\}$.
\end{pf}

\begin{prop}
Let $S$ be a Young diagram. A necessary and sufficient condition for the
existence of a rational point $U\in F_{A_S}$ such that $U\cdot U=U$, is
that
$0\notin S$ and $\Z-S$ should be closed under addition.
\end{prop}

\begin{pf}
Obvious.
\end{pf}

This condition will be called the Weierstrass gap property (WGP). Let us
denote by
$\gr_W(V)$ the open subscheme of $\grv$ consisting of the union of the open
subsets $F_{A_S}$ such that $S$ satisfies WGP. Then one has:

\begin{cor}\label{cor:equationofm}
The subset $\{U\in\gr^n(V)\,\vert\,k\subset U\text{ and } U\cdot
U\subseteq U\}$ is given by one of the following equivalent conditions:
\begin{enumerate}
\item $$\gr_W(V)\cap\left\{ U\in\grv\,\vert\,
\res_{z=0}\psi_U(z,t)\psi_U(z,t')\psi^*_U(z,t'')z^{n-3}dz=0
\,\right\}$$
\item $$\left\{ \begin{aligned}
\res_{z=0}&\psi^*_U(z,t)\frac{dz}{z^{n+1}}\,=\,0
\\  \res_{z=0}&\psi_U(z,t)\psi_U(z,t')\psi^*_U(z,t'')z^{n-3}dz\,=\,0
\end{aligned}\right.$$
\end{enumerate}
\end{cor}

\begin{pf}
The first one is obvious. For the second one we only need to show that
the condition $k\subset U$ is equivalent to
$\res_{z=0}\psi^*_U(z,t)\frac{dz}{z^{n+1}}=0$.
However, from the proof of Theorem
{\ref{thm:res-equa}} and formula \ref{eq:BA=sum} is easily deduced that
$\res_{z=0}f(z)\cdot\psi^*_U(z,t)\frac{dz}{z^{n+1}}=0$ if and only if
$f(z)\in U$.
\end{pf}

\begin{rem}
Assume that $U\in F_{A_S}\subset\gr^n(V)(\spk)$ ($S$ the sequence
associated to a Young diagram) is a point that lies on the image of the
Krichever morphism; that is, there exists
$(C,p,\phi)\in\M_\infty^g$ such that
$K(C,p,\phi)=U$. Note that by the very construction one has an
isomorphism $H^0(C-p,\o_C)\iso U$ and hence:
\begin{itemize}
\item $n=1-g$,
\item the arithmetic genus of $C$ equals $\#({\mathbb Z}_{<0}\cap S)$,
\item the set of Weierstrass gaps of $C$ at $p$ is precisely ${\mathbb
Z}_{<0}\cap S$.
\end{itemize}
\end{rem}


\begin{thm}\label{thm:pde-mod}
The condition:
$$\res_{z=0}\psi_U(z,t)\psi_U(z,t')\psi^*_U(z,t'')
\frac{dz}{z^{2+g}}\,=\,0$$
for a rational point $U\in\gr^{1-g}(V)$ is
equivalent to the infinite set of equations:
$$P(\lambda_1,\lambda_2,\lambda_3)\vert \Sb t=0 \\ t'=0 \\ t''=0 \endSb
\left(\tau_U(t)\cdot \tau_U(t')\cdot\tau_U(t'')\right)\,=\,0$$
for every three Young diagrams $\lambda_1,\lambda_2,\lambda_3$, where
$P(\lambda_1,\lambda_2,\lambda_3)$ is the differential operator defined
by:
$$\sum
p_{\beta_1}(\tilde\partial_t)D_{\lambda_1,\alpha_1}(-\tilde\partial_t)\cdot
p_{\beta_2}(\tilde\partial_{t'})D_{\lambda_2,\alpha_2}(-\tilde\partial_{t'})
\cdot
p_{\beta_3}(-\tilde\partial_{t''})D_{\lambda_3,\alpha_3}(\tilde\partial_{t''})
$$
where the sum is taken over the 6-tuples
$\{\alpha_1,\beta_1,\alpha_2,\beta_2,\alpha_3,\beta_3\}$ of integers such that
$-\alpha_1+\beta_1-\alpha_2+\beta_2-\alpha_3+\beta_3=1+g$.
\end{thm}

\begin{rem}
The meaning of the Residue
Identity $Res_{z=0}\psi\psi\psi^*z^{-(g+2)}dz=0$ (where $\psi$ is the
Baker function of  $(C,p,z)$) is the following: for all sections
$s_1,s_2\in U=H^0(C-p,\o_C)$, $\omega\in
U^\perp=H^0(C-p,\Omega_C)$ the differential $s_1\cdot
s_2\cdot \omega$ has residue zero at $p$. Or, what
amounts to the same: let $D_i$ be the divisor of
poles of $s_i$ ($i=1,2$) and $D^*$ that of $\omega$,
then $D_1+D_2+D^*\,=\, K+m\cdot p$ for some non negative integer $m$
and some canonical divisor $K$.
\end{rem}

These differential equations are the equations of the moduli of
curves (the image of the functor $K$) in the infinite Grassmannian.
A very important fact is that a theta function of a Jacobian variety
satisfies these differential equations, which are not obtained from the KP
equations. Moreover:

\begin{cor}\label{cor:pde-tau}
A formal series $\tau(t)\in k\{\{t_1,t_2,\dots\}\}$ is the
$\tau$-func\-tion associated with a rational point of
$\M_\infty^g\subset \gr^{1-g}(V)$ (and  may therefore be written in terms
of the theta function of a Jacobian variety) if and only if it satisfies
the following set of equations:
\begin{enumerate}
\item the KP equations (given in theorem {\ref{thm:KP-eq}}),
\item the p.d.e.'s  given in theorem {\ref{thm:pde-mod}},
\item the p.d.e.'s:
$$\left(\sum_{-\alpha+\beta=1-g}
p_{\beta}(-\tilde\partial_t)D_{\lambda,\alpha}(\tilde\partial_t)
\right)\vert_{t=0}\tau_U(t)=0\quad
\text{for all Young diagrams }\lambda$$
\end{enumerate}
\end{cor}

\begin{pf}
Note that the third condition is
$\res_{z=0}\psi^*_U(z,t)\frac{dz}{z^{2-g}}=0$ but given in terms of
partial differential equations.
\end{pf}

\begin{rem}
These technics have been used in \cite{Pl} for the study of the moduli
space of Prym varieties and to generalize the characterizations of
Jacobians given by Mulase (\cite{Mul}) and Shiota (\cite{Sh})
\end{rem}

\begin{rem}
An open problem now is to re--state Corollary \ref{cor:equationofm} as
a characterization for a pseudodifferential operator to come from
algebro-geometric data.
\end{rem}



\end{document}